\documentclass[11pt]{article}

\usepackage{amsmath}
\usepackage{amsfonts}
\usepackage{amssymb}
\usepackage{amsthm}
\usepackage{breqn}
\usepackage{setspace}
\usepackage{fullpage}
\usepackage{enumitem}
\usepackage{bbold} 
\usepackage{comment}
\usepackage{hyperref}
\usepackage{bbm}
\usepackage{tikz}
\usepackage{enumitem}
\usepackage{mathtools} 
\usetikzlibrary{patterns,arrows,decorations.pathreplacing}

\newtheorem{observation}{Observation} 
\newtheorem{theorem}{Theorem} 
\newtheorem{lemma}[theorem]{Lemma}
 
\newtheorem{definition}[theorem]{Definition}
\newtheorem{claim}{Claim}

\usepackage{authblk}

\usepackage{authblk}

\begin{document}
\title{The Minimum Number of $4$-Cycles in a Maximal Planar Graph with Small Number of Vertices}

\author[1,2]{Ervin Gy\H{o}ri} 
\author[2]{Addisu Paulos}
\author[2,3]{Oscar Zamora} 

\affil[1]{Alfr\'ed R\'enyi Institute of Mathematics, Hungarian Academy of Sciences. \par
\texttt{gyori.ervin@renyi.mta.hu}}
\affil[2]{Central European University, Budapest.\par
\texttt{ oscarz93@yahoo.es, addisu_2004@yahoo.com}}
\affil[3]{Universidad de Costa Rica, San Jos\'e.}
\maketitle
\begin{abstract}
Hakimi and Schmeichel determined a sharp lower bound for the number of cycles of length 4 in a maximal planar graph with $n$ vertices, $n\geq 5$. It has been shown that the bound is sharp for $n = 5,12$ and  $n\geq 14$ vertices. However, the authors only conjectured the minimum number of cycles of length 4 for maximal planar graphs with the remaining small vertex numbers. In this note, we confirm their conjecture.
\end{abstract}
\section{Introduction}
A \emph{planar graph} is a graph that can be embedded in the plane, i.e., it can be drawn on the plane in such a way that no two edges intersect except at their end vertices. Such a drawing is a \emph{plane graph}. A \emph{maximal planar} is a planar graph such that adding an edge between any two non-adjacent vertices results in a non-planar graph. 

For a maximal planar graph $G$, denote $C_4(G)$ to be the number of 4-cycles in the graph. Hakimi and Schmeichel \cite{hakimi} obtained a lower bound for $C_4(G)$ when the number of vertices is $n\geq 5$.
It has been shown that for any maximal planar graph $G$ on $n$ vertices, with $n=5, 12$ or $n\geq 14$, $C_4(G)\geq 3n-6$. 
Indeed, if $n=5$, there is only a unique maximal planar graph and the number of $4$-cycles is 9. 
For $n\geq 12$ and $n\neq 13$, it has been shown in \cite{hakimi2} that there is a 5-connected maximal planar graph with the indicated number of vertices. 
Since after deleting an edge in a maximal planar graph we obtain a face bounded by a 4-cycle, and each edge represents a unique $4$-face (the only case where two edges represent the same face is when $G=K_4$), and so $C_4(G)\geq 3n-6$. Now if there is any other cycle of length 4 in $G$, then it should be a separating $4$-cycle, i.e., a cycle of length three such that there are vertices of $G$ in both the exterior and interior region of the cycle. However, there is no 4-vertex cut set as the graph is 5-connected. Therefore, the bound is sharp.
 
In \cite{hakimi} the authors indicated their conjecture about the minimum number of $4$-cycles in the maximal planar graphs of $n$ vertices, where $n=4$ or $6\leq n\leq 11$ or $n=13$. In this note, we confirm their conjecture for each such value of $n$. In addition to the conjecture, the authors gave the maximal planar graphs attaining each value indicated in their conjecture.

The following notations and terminologies are important. For a graph $G$, $V(G)$ and $E(G)$ denote the vertex and edge sets of $G$ respectively. $v(G)$ and $e(G)$ denote the number of vertex and  and edge in $G$ respectively. For a vertex $v$ in $G$, $d_G(v)$ denotes the degree of $v$, and $N(v)$ denotes the set of vertices adjacent to $v$. We may omit the subscript if the underlying graph is obvious. $\delta(G)$ and $\Delta(G)$ denote the minimum and maximum degrees in $G$ respectively. We use the term $k$-cycle to describe a cycle of length $k$. By a $C_k(G)$ we mean the number of $k$-cycles in $G$. For $v\in V(G)$, $C_k(G,v)$ denote the number of $4$-cycles containing $v$. In a plane graph $G$, by \emph{separating $k$-cycle} in $G$, we mean a cycle of length $k$ such that there are vertices of $G$ in both the exterior and interior region of the cycle. We denote a $k$-cycle with vertices $v_1, \ v_2,\ \dots,\ v_k$ in sequential order by $(v_1,\ v_2,\ \dots\, v_k,\ v_1)$. A cycle $\mathcal{C}$ in $G$ is chordal if $uw\in E(G)$ for some $u,\ w\in V(\mathcal{C})$ and $uw\notin E(\mathcal{C})$. For a positive integer $k$, $[k]$ is the set of the first $k$ positive integers. 

It should be noticed that, every face in a plane drawing of a maximal planar graph is a $3$-face. If $G$ is a maximal planar graph, $e(G)=3n-6$. Moreover, for a vertex $v\in V(G)$ with $d(v)=k$, then $N(v)$ induces a unique $k$-cycle. If $v(G)\geq 4$ we have, $3\leq \delta(G)\leq 5$. If $\delta(G)=5$, then $v(G)\geq 12$, i.e., the minimum number of vertices required for a maximal planar graph with a minimum degree $5$ is at least $12$.

\begin{definition}
Define $g(n,\ C_4)$ as the minimum number of 4-cycles in a maximal planar graph of $n$ vertices,
i.e., $$g(n,C_4):=\min\left\{\ C_4(G)\ |\ \text{ $G$ is an $n$-vertex maximal planar graph  }\right\}.$$
\end{definition}

From the result of Hakimi and Schmeichel, for $n=5,\ n\geq 12$ and $n\neq 13,\ g(n,C_4)=3n-6$. 
In this note, we determined the exact value of $g(n,C_4)$ for each $n\in \{4,\ 6,\ 8,\ 9,\ 10,\ 11,\ 13 \}$. 
The following theorem gives the values.

\begin{theorem}\ \\
\label{thm1}
\begin{align*}
g(n,C_4)= \begin{cases}
3, &\text{if $n=4$;}\\
15, &\text{if $n=6$;}\\
20, &\text{if $n=7$;}\\
23, &\text{if $n=8$;}\\
24, &\text{if $n=9$;}\\
26, &\text{if $n=10$;}\\
29, &\text{if $n=11$;}\\
34, &\text{if $n=13$;}\\
3n-6, &\text{if otherwise.}
\end{cases}
\end{align*}
\end{theorem}
We need the following three lemmas to prove our main results. The proof of Lemma~\ref{ind} and Lemma~\ref{ind1} relies on the property of maximal outerplanar graphs as stated in Lemma~\ref{ind2}.
\begin{lemma}\cite{west}\label{ind2}
Every maximal outerplanar graph with at least $3$ vertices contains a degree-$2$ vertex.
\end{lemma}
\begin{lemma}\label{ind}
Let $G$ be an $n$-vertex maximal plane graph, where $n\geq 4$. If $\Delta(G)=n-1$, then there is a vertex of degree 3 in $G$. 
\end{lemma}
\begin{proof}
Let $x\in V(G)$ such that $d(x)=n-1$. Since $G-x$ is a maximal outerplanar graph, then there is a degree-$2$ vertex, say $v$, in $G-x$. Thus, $d(v)=3$.
\end{proof}
\begin{lemma}\label{ind1}
Let $G$ be an $n$-vertex maximal plane graph, where $n\geq6$. Let $v\in V(G)$ such that $d(v)=n-2$. If the cycle $\mathcal{C}$ induced by $N(v)$ has a chord, then there is a vertex of degree 3 in $G$.
\end{lemma}
\begin{proof}
Let  $e=\{x,y\}\in E(G)$ be a chord in $\mathcal{C}$, where $x,\ y\in N(v)$. The edge $e$ cuts $\mathcal{C}$ into two non-overlapping paths, $P_1$ and $P_2$ with common end vertices. Let $R_1$ be the region enclosed by $P_1$ and $e$ and not containing $v$, and $R_2$ be the region enclosed by $P_2$ and $e$ and not containing $v$. 

Denote $u\in V(G)$ such that $u\notin N(v)$.  Without loss of generality, we may assume that $u$ is in $R_1$. Deleting $v$, $u$, and all the interior vertices of $P_1$, then the remaining graph is a maximal outerplanar graph, say $G'$. Hence we have a vertex $w\in V(G')$ such that $d_{G'}(w)=2$. Notice that $w\in N(v)$ and an interior vertex of $P_2$. Thus, $d_G(w)=3$.     
\end{proof}
\begin{lemma}\label{lema}
Let $G$ be an $n$-vertex maximal plane, where $n\geq 7$. If $d_G(v_1)=d_G(v_2)=4$ for some $v_1,\ v_2 \in V(G)$, then $N(v_1)$ and $N(v_2)$ induce two different $4$-cycles.
\end{lemma}
\begin{proof}
Suppose for contradiction that $N(v_1)$ and $N(v_2)$ induces the same 4-cycles. 
Then one vertex, say $v_1$,  is in the interior the 4-cycle and joined to all the four vertices of the cycle. 
Hence the other vertex, $v_2$ is on the exterior region of the cycle and is incident to every vertex of the cycle.
Since the number of vertices is at least 7, then at least one vertex is not used. If there is a vertex $u\not \in N(v_1)\cup N(v_2)\cup \{v_1,\ v_2\}$ in the interior region of the 4-cycle, then there must be another vertex other than the vertices in $N(v_1)$ that is incident to $v_1$, which is a contradiction to the fact that $d_G(v_1)=4$. 
\end{proof}

\section{Proof of Theorem \ref{thm1}}

Since there is only one maximal planar graph with 4 vertices, it is easy to see that the number of $4$-cycles in the graph is exactly 3, and hence $g(4,C_4)=3$. 

On the other hand, it can be checked that there are only two $6$-vertex maximal plane graphs.  Indeed, let $G$ be a $6$-vertex maximal plane graph. First, we see that there is a separating $4$-cycle in $G$. Observe that if there is a vertex of degree 3, say $v$, then there are two vertices outside of the $3$-cycle induced by $N(v)$. Thus, we obtain separating $4$-cycle taking the vertices of the cycle and a vertex incident to the face which uses an edge of the cycle and lying in the region where the two vertices are. If there is no vertex of degree 3, then we have a $w\in V(G)$ such that $d_G(w)=4$, and obviously the cycle induced by $N(w)$ forms a separating $4$-cycle. 

Let $(v_1,\ v_2,\ v_3,\ v_4)$ be a separating $4$-cycle in $G$. Let the interior and exterior vertices with reference to the cycle be $u$ and $v$ respectively. Then due to symmetry property, we consider when degree pairs of $u$ and $v$ are either $(4,\ 4)$ or $(4,\ 3)$ or $(3,\ 3)$ respectively. 

For the first degree pair, $N(u)=N(v)=\{v_1,\ v_2,\ v_3,\ v_4\}$. In this case, the maximal planar graph is shown in Figure~\ref{10} (left), and $C_4(G)=15$.

For the second degree pair, observe that $N(u)=\{v_1,\ v_2,\ v_3,\ v_4\}$. Without loss of generality, we may assume that $N(v)=\{v_1,\ v_2,\ v_3\}.$ Since $d_G(v)=3$, the $v_1v_3\in E(G)$. In this case, the resulting graph is isomorphic to the maximal planar graph shown in Figure~\ref{10} (right) and $C_4(G)=16$. 

Finally, we assume the third degree pair, i.e., $d_G(u)=d_G(v)=3$. Let $N(u)=\{v_1,\ v_2,\ v_3\}$. This implies, $v_1v_3\in E(G)$ and $v_1v_3$ is incident to the $3$-face bounded by the path $(v_1,\ u,\ v_3)$. Since $\delta(G)\geq 3$, then $v_4\in N(v).$ If both $v_1,\ v_3\in N(v)$, then $v_1v_3$ must be incident to the $3$-face bounded by the path $(v_1,\ v,\ v_3)$. But this results in a contradiction as the graph is not a multi-graph. Thus, without loss of generality, we may assume that $N(v)=\{v_4,\ v_3,\ v_2\}$. In this case, $v_2v_4\in  E(G)$ and $v_2v_4$ is incident to the $3$-face bounded by the path $(v_2,\ v,\ v_4)$. Therefore, $d_G(v_4)=4$ and hence we can take the separating $4$-cycle induced by $N(v_4)$. But this is already a settled case. 

Therefore, there are only two maximal plane graphs with 6 vertices and $g(6,C_4)=15$.

\begin{figure}[h]
\centering
\begin{tikzpicture}[scale=0.25]
\draw[ultra thick](0,10)--(8.7,-5)--(-8.7,-5)--(0,10);
\draw[ultra thick](0,-2.5)--(2.16,1.25)--(-2.16,1.25)--(0,-2.5)(-8.7,-5)--(-2.16,1.25)--(0,10)(8.7,-5)--(2.16,1.25)--(0,10)(-8.7,-5)--(0,-2.5)--(8.7,-5);
\draw[fill=black](0,-2.5)circle(15pt);
\draw[fill=black](2.16,1.25)circle(15pt);
\draw[fill=black](-2.16,1.25)circle(15pt);
\draw[fill=black](0,10)circle(15pt);
\draw[fill=black](8.7,-5)circle(15pt);
\draw[fill=black](-8.7,-5)circle(15pt);
\node at (0,-9) {$C_4(G)=15$};
\end{tikzpicture}\qquad\qquad\qquad
\begin{tikzpicture}[scale=0.25]
\draw[ultra thick](0,10)--(8.7,-5)--(-8.7,-5)--(0,10)(0,10)--(0,-1.25)(-8.7,-5)--(0,6.25)--(8.7,-5)(-8.7,-5)--(0,2.5)--(8.7,-5)(-8.7,-5)--(0,-1.25)--(8.7,-5);
\draw[fill=black](0,10)circle(15pt);
\draw[fill=black](0,6.25)circle(15pt);
\draw[fill=black](0,2.5)circle(15pt);
\draw[fill=black](0,-1.25)circle(15pt);
\draw[fill=black](8.7,-5)circle(15pt);
\draw[fill=black](-8.7,-5)circle(15pt);
\node at (0,-9) {$C_4(G)=16$};
\end{tikzpicture}
\caption{The two $6$-vertex maximal planar graphs with their respective number of $4$-cycles.}
\label{10}
\end{figure}

Next, we prove a sequence of claims to prove the remaining cases of the theorem. The following observation is important in proving the claims.
\begin{observation}\label{ob11}
\emph{Let $G$ be an $n$-vertex maximal planar graph, where $7\leq n\leq 12$ or $n=13.$ If $G$ contains a degree-$3$ vertex, say $v$, then $G-v$ is an $(n-1)$-vertex maximal planar graph. Since $n\geq 5$, it can be checked that $C_4(G,v)\geq 6$, and hence $C_4(G)\geq g(n-1, C_4)+C_4(G,v)\geq g(n-1, C_4)+6.$} 
\end{observation}

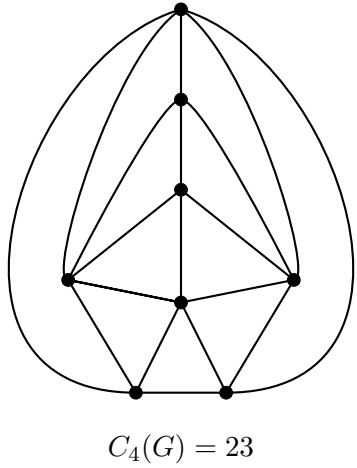
\begin{figure}[h]
\centering
\begin{tikzpicture}[scale=0.25]
\draw[ultra thick](0,10)--(8.7,-5)--(-8.7,-5)--(0,10);
\draw[ultra thick](0,-2.5)--(2.16,1.25)(-2.16,1.25)--(0,-2.5)(-8.7,-5)--(-2.16,1.25)--(0,10)(8.7,-5)--(2.16,1.25)--(0,10)(-8.7,-5)--(0,-2.5)--(8.7,-5)(0,3)--(2.16,1.25)--(0,0.15)--(-2.16,1.25)--(0,3)(0,10)--(0,0.15)--(0,-2.5);
\draw[fill=black](0,3)circle(13pt);
\draw[fill=black](0,0.15)circle(13pt);
\draw[fill=black](0,-2.5)circle(13pt);
\draw[fill=black](2.16,1.25)circle(13pt);
\draw[fill=black](-2.16,1.25)circle(13pt);
\draw[fill=black](0,10)circle(13pt);
\draw[fill=black](8.7,-5)circle(13pt);
\draw[fill=black](-8.7,-5)circle(13pt);
\node at (0,-9) {$C_4(G)=23$};
\end{tikzpicture}
\caption{Maximal planar graph with 8 vertices containing a degree 5 vertex and no vertex of degree 3.}
\label{7}
\end{figure}

\begin{claim}$g(7,C_4)=20.$ \label{cl2}\end{claim}
\begin{proof}
Let $G$ be a $7$-vertex maximal planar graph. If there is a vertex of degree 3, say $v$, then by Observation~\ref{ob11}  $C_4(G)=C_4(G-v)+C_4(G,v)\geq g(6,C_4)+6=21$. It can be checked that $G$ is unique if we assume $\delta(G)=4$. To see this, let $v\in V(G)$ such that $d(v)=4$. Let $(v_1,\ v_2,\ v_3, v_4)$ be the separating $4$-cycle induced by vertices in $N(v)$. Let $u$ and $w$ be the vertices in the other region of the separating cycle not containing $v$. Since $\delta(G)=4$, then both $u$ and $w$ must be adjacent to three vertices in $N(v)$. Without loss of generality, let $N(u)=\{v_1,\ v_2,\ v_3\}$ and $N(w)=\{v_1,\ v_4,\ v_3\}$. For the same reason, $\delta(G)=4$, we have $uw\in E(G)$. Therefore, the graph is uniquely determined. Figure~\ref{9} shows the graph and it has been mentioned in~\cite{hakimi} that the number of $4$-cycles in this graph is $20$. Therefore, $g(7, C_4)=20.$    
\end{proof}
\begin{figure}[h]
\centering
\begin{tikzpicture}[scale=0.25]
\draw[ultra thick](0,10)--(8.7,-5)--(-8.7,-5)--(0,10);
\draw[ultra thick](0,-2.5)--(2.16,1.25)(-2.16,1.25)--(0,-2.5)(-8.7,-5)--(-2.16,1.25)--(0,10)(8.7,-5)--(2.16,1.25)--(0,10)(-8.7,-5)--(0,-2.5)--(8.7,-5)(2.16,1.25)--(0,2.5)--(-2.16,1.25)(0,10)--(0,2.5)--(0,-2.5);
\draw[fill=black](0,2.5)circle(13pt);
\draw[fill=black](0,-2.5)circle(13pt);
\draw[fill=black](2.16,1.25)circle(13pt);
\draw[fill=black](-2.16,1.25)circle(13pt);
\draw[fill=black](0,10)circle(13pt);
\draw[fill=black](8.7,-5)circle(13pt);
\draw[fill=black](-8.7,-5)circle(13pt);
\node at (0,-9) {$C_4(G)=20$};
\end{tikzpicture}
\caption{The only maximal planar graph on 7 vertices with minimum degree 4.}
\label{9}
\end{figure}
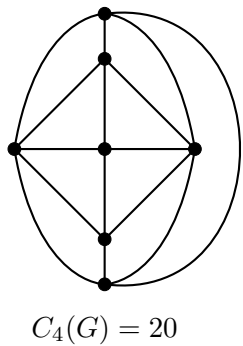
\begin{claim}$g(8,C_4)=23.$ \label{cll3}\end{claim}
\begin{proof}
Let $G$ be an $8$-vertex maximal plane graph. If there is a vertex of degree 3, then by Observation~\ref{ob11} and Claim~\ref{cl2} we have, $C_4(G)\geq g(7,C_4)+6\geq 26$.

Assume $\delta(G)\geq 4$. If there is a degree-$7$ vertex in $G$, then by Lemma \ref{ind}, there is a vertex of degree 3 in $G$ as well, which is a contradiction. Let $G$ have a vertex, say $v$, of degree 6. Then we have a vertex, say $u$, which is inside a region bounded by the induced 6-cycle of $N(v)$, and not containing $v$. If the induced cycle is chordal, then by Lemma \ref{ind1}, there is a degree-$3$ vertex, which is a contradiction to our assumption. If the cycle is not chordal, then $u$ should be adjacent to all of the vertices of the cycle. Thus, all of the vertices of the cycle are of degree 4. Therefore, using Lemma \ref{lema}, we have at least 6 separating 4-cycles. Therefore, $C_4(G)\geq (3n-6)+6=24$ where $n=8$.

Finally we may assume that $4\leq d_G(v)\leq 5$, for all $v\in V(G)$. Let $v$ be a vertex of degree $5$ in $G$. Then there are 2 vertices inside the region (not containing $v$) bounded by the $5$-cycle induced by $N(v)$. Each of the two vertices should be adjacent to at least three vertices of the cycle. Otherwise, one of the two vertices is with degree 3. Moreover, the two vertices are adjacent, and only one maximal plane graph meets the properties, see Figure \ref{7}. Notice that the graph is mentioned in \cite{hakimi}, and the number of $4$-cycles in the graph is $23$. Therefore, we have, $g(8,C_4)=23$.    
\end{proof}

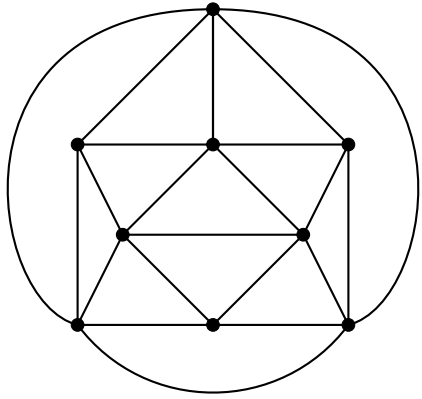
\begin{figure}[h]
\centering
\begin{tikzpicture}[scale=0.25]
\draw[ultra thick](0,10)--(8.7,-5)--(-8.7,-5)--(0,10);
\draw[ultra thick](0,-2.5)(2.16,1.25)(-2.16,1.25)(0,-2.5)(-8.7,-5)--(-2.16,1.25)--(0,10)(8.7,-5)--(2.16,1.25)--(0,10)(-8.7,-5)--(0,-2.5)--(8.7,-5)(0,3)--(2.16,1.25)(-2.16,1.25)--(0,3)--(0,10)(2.16,1.25)--(2.16,-1)--(0,-2.5)--(-2.16,-1)--(-2.16,1.25)(8.7,-5)--(2.16,-1)(-8.7,-5)--(-2.16,-1)--(2.16,-1)--(0,3)--(-2.16,-1);
\draw[fill=black](2.16,-1)circle(13pt);
\draw[fill=black](-2.16,-1)circle(13pt);
\draw[fill=black](0,3)circle(13pt);
\draw[fill=black](0,-2.5)circle(13pt);
\draw[fill=black](2.16,1.25)circle(13pt);
\draw[fill=black](-2.16,1.25)circle(13pt);
\draw[fill=black](0,10)circle(13pt);
\draw[fill=black](8.7,-5)circle(13pt);
\draw[fill=black](-8.7,-5)circle(13pt);
\node at (0,-9) {$C_4(G)=24$};
\end{tikzpicture}
\caption{Maximal planar graph with 9 vertices and with the least number of 4-cycles.}
\label{5}
\end{figure}

\begin{claim}$g(9,C_4)=24.$ \label{cl3}\end{claim}
\begin{proof}
Let $G$ be a $9$-vertex maximal plane graph. If there is a vertex of degree 3, then by Observation~\ref{ob11} and Claim~\ref{cll3} we have $C_4(G)\geq g(8,\ C_4)+6=29.$

Assume $\delta(G)\geq 4$. If there is a degree-$8$, then, by Lemma \ref{ind}, we have a vertex of degree 3 in the graph, which is a contradiction to our assumption. Let $G$ have a $v$ such that $d_G(v)=7$. Then there is a vertex, say $u$, in the region bounded by the $7$-cycle induced by $N(v)$ and not containing $v$. If the cycle is chordal, then, by Lemma \ref{ind1}, we have a vertex of degree $3$ in the graph, which is impossible as $\delta(G)\geq 4$ by assumption. Therefore, $u$ must be incident to all vertices in $N(v)$. Since there are 7 degree-$4$ vertices, then we have at least 7 separating 4-cycles. Thus, $C_4(G)\geq (3n-6)+7=28$.

Now assume $4\leq d_G(v)\leq 6$, for all $v\in V(G)$. We claim that $G$ has at least 3 vertices of degree-$4$. Indeed, let $n_4,\ n_5$ and $n_6$ be the number vertices of degree $4$, $5$ and $6$ respectively. Then $n_4+n_5+n_6=9$. Moreover, from the property that $\sum\limits_{v\in V(G)}d(v)=2e(G)$ and $e(G) =3n-6$, we have $4n_4+5n_5+6n_6=42$. So, solving the two equations simultaneously, we get $n_6=n_4-3$. Therefore, $G$ contains at least 3 vertices of degree 4. In other words, $G$ has at least 3 distinct separating 4-cycles. Thus, $C_4(G)\geq (3n-6)+3=24$. From the example given in Figure~\ref{5}, which is also mentioned in \cite{hakimi}, and considering all the bounds mentioned above we have, $g(9,C_4)=24$.    
\end{proof}

\begin{figure}[h]
\centering
\begin{tikzpicture}[scale=0.25]
\draw[ultra thick](0,10)--(8.7,-5)--(-8.7,-5)--(0,10);
\draw[ultra thick](0,-2.5)(2.16,1.25)(-2.16,1.25)(0,-2.5)(-8.7,-5)--(-2.16,1.25)--(0,10)(8.7,-5)--(2.16,1.25)--(0,10)(-8.7,-5)--(0,-2.5)--(8.7,-5)(0,3)--(2.16,1.25)(-2.16,1.25)--(0,3)--(0,10)(2.16,1.25)--(2.16,-1)--(0,-2.5)--(-2.16,-1)--(-2.16,1.25)(8.7,-5)--(2.16,-1)(0,0.5)--(0,3)--(2.16,1.25)(0,3)--(-2.16,1.25)(2.16,-1)--(-2.16,-1)(0,0.5)--(2.16,-1)(0,0.5)--(-2.16,-1)(0,0.5)--(2.16,1.25)(0,0.5)--(-2.16,1.25)(-8.7,-5)--(-2.16,-1);
\draw[fill=black](0,0.5)circle(13pt);
\draw[fill=black](2.16,-1)circle(13pt);
\draw[fill=black](-2.16,-1)circle(13pt);
\draw[fill=black](0,3)circle(13pt);
\draw[fill=black](0,-2.5)circle(13pt);
\draw[fill=black](2.16,1.25)circle(13pt);
\draw[fill=black](-2.16,1.25)circle(13pt);
\draw[fill=black](0,10)circle(13pt);
\draw[fill=black](8.7,-5)circle(13pt);
\draw[fill=black](-8.7,-5)circle(13pt);
\node at (0,-9) {$C_4(G)=26$};
\end{tikzpicture}
\caption{Maximal planar graph with 10 vertices and with the least number of 4-cycles.}
\label{4}
\end{figure}
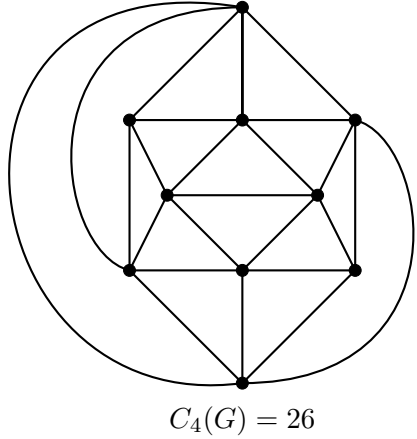

\begin{claim}$g(10,C_4)=26.$ \end{claim}
\begin{proof}
Let $G$ be a maximal plane graph with $10$ vertices. If there is a vertex of degree 3, then by Observation~\ref{ob11} and Claim~\ref{cl3} we have $C_4(G)\geq g(9, C_4)+6\geq 30$.

Suppose that $\delta(G)\geq4$. If there is a vertex of degree 9, then by Lemma \ref{ind}, $G$ contains a degree 3 vertex. If there is a vertex of degree 8, say $v$, then considering the $8$-cycle induced by $N(v)$, we have a vertex, say $u$, inside of the region induced by the cycle and not containing $v$. If the cycle is chordal, then by Lemma \ref{ind1}, there is a vertex of degree 3 in the graph. Otherwise, $u$ is adjacent to every vertex in $N(v)$. This implies, that each of the 8 vertices of the cycle is of degree 4. Thus, $G$ has at least 8 separating 4-cycles. Therefore, $C_4(G)\geq (3n-6)+8=32$, where $n=10$. If there is a vertex of degree 7, say $v$,  then considering the 7-cycle induced by $N(v)$, there will be two vertices that lie in the region bounded by the 7-cycle not containing $v$. Each of these vertices should be adjacent to at least $3$ vertices of the cycle. Otherwise, there exists a vertex of degree 3. Thus, there exist at least two vertices of degree 4 on the graph. This implies we have at least 2 separating 4-cycles. Thus, $C_4(G)\geq (3n-6)+2\geq 26$, where $n=10$.

Now for the remaining cases we assume $4\leq d_G(v)\leq 6$,  
 for all $v\in V(G)$,  we show that there are at least 2 vertices of degree 4. Let $n_4$, $n_5$, and $n_6$ be the number of vertices of degree 4, 5, and 6 respectively. Then $n_4+n_5+n_6=10$ and from the property that $\sum\limits_{v\in V(G)}d(v)=2e(G)$ and $e(G)=3n-6$ and $n=10$, we have $4n_4+5n_5+6n_6=48$. Equating the two equations simultaneously we get, $n_6=n_4-2$. 
Thus, the graph contains at least two degree 4 vertices. Therefore, $C_4(G)\geq 26$. This value is sharp, due to the graph in Figure \ref{4} (also mentioned in \cite{hakimi}). Considering all the results we obtained, we have $g(10,C_4)=26$.    
\end{proof}

\begin{claim}$g(11,C_4)=29.$ \end{claim}
\begin{proof}
Let $G$ be an $11$-vertex maximal plane graph. If there is a vertex of degree 3 in the graph, then from the fact that $g(10,C_4)=26$ we have, $C_4(G)\geq g(10,C_4)+6\geq 32$.

Assume that $\delta(G)\geq 4$. If there is a vertex of degree $10$, then by Lemma~\ref{ind} the graph contains a vertex of degree 3. If there is a vertex of degree 9 in the graph, then following a similar argument given as before, it can be checked that either there is a vertex of degree 3 or $C_4(G)\geq 36$. If there is a vertex, say $v$, of degree 8, then there are two vertices that lie in a region bounded by the 8-cycle induced by $N(v)$ and not containing $v$. The vertices are adjacent to at least 3 vertices of the cycle. Otherwise, one of the two vertices becomes a degree-$3$ vertex, which is not true by our assumption $\delta(G)\geq 4$. Since the vertices are adjacent to at least 3 vertices of the cycle, then there are two vertices of degree 4 or there is a vertex of degree 3. Thus, we have at least two separating $4$-cycles. Therefore, $C_4(G)\geq 29$.

Now we may assume that $4\leq d(v)\leq 7$, for all $v\in V(G)$. Let $n_4,\ n_5,\ n_6$, and $n_7$ be the number of vertices of the graph with degrees 4, 5, 6 and 7 respectively. Since the number of vertices is 11, then $n_4+n_5+n_6+n_7=11$. Moreover, from the property that $\sum\limits_{v\in V(G)}d(v)=2e(G)$ and $e(G)=3n-6$ and $n=11$, we have $4n_4+5n_5+6n_6+7n_7=54$. Solving the two equations simultaneously we get, 
\begin{align}\label{jnk}
  n_6+2n_7=n_4-1  
\end{align}

If there is a vertex of degree $7$ in $G$, then $n_7>0$. Since $n_6\geq 0$ we have $n_4\geq 3$. That means, there are at least 3 separating 4-cycles. Thus, $C_4(G)\geq 30$. 

If there is no vertex of degree 7, then $4\leq d(v)\leq 6$ for all $v\in V(G)$ and $n_7=0$. In this case, we claim that there are at least 2 degree 4 vertices. To see this, suppose that there is at most one degree 4 vertex. The remaining vertices are of degrees 5 and 6. From equation~\ref{jnk}, we have $n_6=n_4-1$. If $n_4=0$, it results in a contradiction. If $n_4=1$, then $n_6=0$ and the remaining vertices are of degree $5$.  Let $e=vy\in E(G)$, where $d_G(v)=4$. Let two vertices $x_1,\ x_2\in N(v)\cap N(y)$ be such that $(v,\ x_1,\ y,\ v)$ and $(v,\ x_2,\ y,\ v)$ are faces. Let the remaining vertex in $G$ which is incident to $v$ be $z$. Here both $x_1$ and $x_2$ are incident to $z$. Otherwise, the degree of $v$ is more than 4.  Consider the vertices in $G$ incident to $y$ be $y_1,\ y_2$ such that $N(y)=\{v,\ x_1,\ x_2,\ y_1,\ y_2\}$ and the 5-cycle induced by $N(y)$ is $(v,\ x_1,\ y_1,\ y_2,\ x_2,\ v)$  as shown on the Figure \ref{3}.

\begin{figure}[ht]
\centering
\begin{tikzpicture}[scale=0.35]
\draw[ultra thick](0,0)--(4,0)--(6,4)(4,0)--(2,-3)(0,0)--(2,3)--(-4,0)--(0,0)--(2,-3)(-4,0)--(-6,4)(-4,0)--(-6,-4)(-4,0)--(2,-3)(2,3)--(4,0)--(6,-4)(-6,4)--(2,3)--(6,4)(-6,-4)--(2,-3)--(6,-4);
\draw[black,ultra thick](6,4)..controls (7,3) and (7,-3) .. (6,-4);
\draw[black,ultra thick](-6,4)..controls (-7,3) and (-7,-3) .. (-6,-4);
\draw[black,ultra thick](-6,4)..controls (-3,6) and (3,6) .. (6,4);
\draw[black,ultra thick](-6,-4)..controls (-3,-6) and (3,-6) .. (6,-4);
\draw[fill=black](0,0)circle(10pt);
\draw[fill=black](2,3)circle(10pt);
\draw[fill=black](2,-3)circle(10pt);
\draw[fill=black](4,0)circle(10pt);
\draw[fill=black](6,4)circle(10pt);
\draw[fill=black](6,-4)circle(10pt);
\draw[fill=black](-4,0)circle(10pt);
\draw[fill=black](-6,-4)circle(10pt);
\draw[fill=black](-6,4)circle(10pt);
\node at (0.8,0.5) {$v$};
\node at (7,-4) {$y_1$};
\node at (7,4) {$y_2$};
\node at (-7,-4) {$z_1$};
\node at (-7,4) {$z_2$};
\node at (-4,-1) {$z$};
\node at (4,-1) {$y$};
\node at (2,4) {$x_2$};
\node at (2,-4) {$x_1$};
\end{tikzpicture}
\caption{Constructing a maximal planar graph on 11 vertices and with only one vertex of degree 4 and all the remaining vertices of degree 5.}
\label{3}
\end{figure}
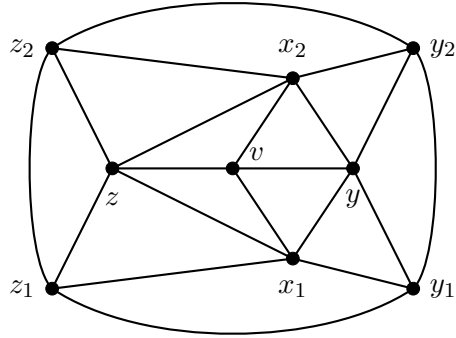
We claim that both $y_1$ and $y_2$ are not incident to the vertex $z$. Indeed, suppose $z$ is adjacent to $y_1$. Consider the triangular region bounded $(z,\ x_1,\ y_1,\ z)$ not containing $v$. If this region contains no vertex, then $d(x_1)=4$, which is a contradiction as $v$ is the only vertex in $G$ with degree 4. If the region contains one vertex, then that vertex is of degree 3, which is again a contradiction. If it contains more than two vertices, then one of the three vertices, $z,\ x_1$ or $y_1$, will be of degree at least 6, which is again a contradiction. 

Let $z_1$ and $z_2$ be the two vertices in $G$ adjacent to $z$ such that $N(z)=\{v,\ x_1,\ x_2,\ z_1,\ z_2\}$ and the 5-cycle induced by $N(z)$ is $(v,\ x_1,\ z_1,\ z_2,\ x_2,\ v)$ as shown in Figure \ref{3}. Since $d(x_1)=d(x_2)=5$, then $z_1y_1$ and $z_2y_2$ are in $E(G)$. Since $(y_1,\ z_1,\ z_2,\ y_2,\ y_1)$ is a separating cycle containing two vertices outside and their degree must be 5, each of them should be adjacent to all the four vertices of the cycle. This implies degrees of the vertices on the cycle are more than 5, which is a contradiction.

Therefore, there are at least 2 vertices of degree 4. Hence $C_4(G)\geq 29$. So with the graph given in Figure~\ref{2}(left), which is also mentioned in the paper \cite{hakimi}), we conclude $g(11, C_4)=29$.
\end{proof}
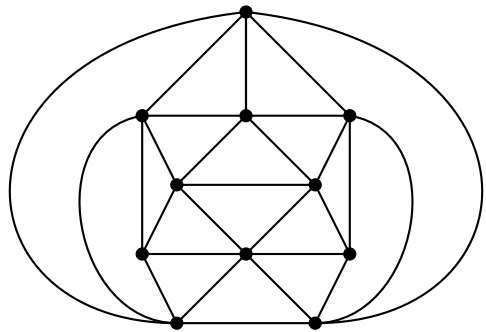
\begin{figure}[h]
\centering
\begin{tikzpicture}[scale=0.25]
\draw[ultra thick](0,10)--(8.7,-5)--(-8.7,-5)--(0,10);
\draw[ultra thick](0,-2.5)(2.16,1.25)(-2.16,1.25)(0,-2.5)(-8.7,-5)--(-2.16,1.25)--(0,10)(8.7,-5)--(2.16,1.25)--(0,10)(-8.7,-5)--(0,-2.5)--(8.7,-5)(0,3)--(2.16,1.25)(-2.16,1.25)--(0,3)--(0,10)(2.16,1.25)--(2.16,-1)--(0,-2.5)--(-2.16,-1)--(-2.16,1.25)(8.7,-5)--(2.16,-1)(0,3)--(2.16,1.25)(0,3)--(-2.16,1.25)(0.9,0)--(0,3)(0.9,0)--(2.16,1.25)(0.9,0)--(2.16,-1)(0.9,0)--(0,-2.5)(-0.9,0)--(0,3)(-0.9,0)--(-2.16,1.25)(-0.9,0)--(-2.16,-1)(-0.9,0)--(0,-2.5)(0.9,0)--(-0.9,0)(-8.7,-5)--(-2.16,-1);
\draw[fill=black](-0.9,0)circle(13pt);
\draw[fill=black](0.9,0)circle(13pt);
\draw[fill=black](2.16,-1)circle(13pt);
\draw[fill=black](-2.16,-1)circle(13pt);
\draw[fill=black](0,3)circle(13pt);
\draw[fill=black](0,-2.5)circle(13pt);
\draw[fill=black](2.16,1.25)circle(13pt);
\draw[fill=black](-2.16,1.25)circle(13pt);
\draw[fill=black](0,10)circle(13pt);
\draw[fill=black](8.7,-5)circle(13pt);
\draw[fill=black](-8.7,-5)circle(13pt);
\node at (0,-9) {$C_4(G)=29$};
\end{tikzpicture}\qquad\qquad
\begin{tikzpicture}[scale=0.25]
\draw[ultra thick](0,10)--(8.7,-5)--(-8.7,-5)--(0,10);
\draw[ultra thick](0,-2.5)(2.16,1.25)(-2.16,1.25)(0,-2.5)(-8.7,-5)--(-2.16,1.25)--(0,10)(8.7,-5)--(2.16,1.25)--(0,10)(-8.7,-5)--(0,-2.5)--(8.7,-5)(0,3)--(2.16,1.25)(-2.16,1.25)--(0,3)--(0,10)(2.16,1.25)--(2.16,-1)--(0,-2.5)--(-2.16,-1)--(-2.16,1.25)(8.7,-5)--(2.16,-1);
\draw[ultra thick](2.16,-1)--(0,-1.3)--(-2.16,-1)(2.16,-1)--(1,1)--(0,3)(-2.16,-1)--(-1,1)--(0,3)(1,1)--(2.16,1.25)(-1,1)--(-2.16,1.25)(0,-2.5)--(0,-1.3)--(1,1)--(-0.3,0)--(0,3)(-0.3,0)--(-1,1)(-0.3,0)--(-2.16,-1)(-0.3,0)--(0,-1.3)(-2.16,-1)--(-8.7,-5);
\draw[fill=black](-1,1)circle(10pt);
\draw[fill=black](1,1)circle(10pt);
\draw[fill=black](0,-1.3)circle(10pt);
\draw[fill=black](-0.3,0)circle(10pt);
\draw[fill=black](2.16,-1)circle(10pt);
\draw[fill=black](-2.16,-1)circle(10pt);
\draw[fill=black](0,3)circle(10pt);
\draw[fill=black](0,-2.5)circle(10pt);
\draw[fill=black](2.16,1.25)circle(10pt);
\draw[fill=black](-2.16,1.25)circle(10pt);
\draw[fill=black](0,10)circle(10pt);
\draw[fill=black](8.7,-5)circle(10pt);
\draw[fill=black](-8.7,-5)circle(10pt);
\node at (0,-9) {$C_4(G)=34$};
\end{tikzpicture}
\caption{Maximal planar graph with 11 and $13$ vertices and with the least number of 4-cycles.}
\label{2}
\end{figure}

\begin{claim}$g(13,C_4)=34.$ \end{claim}
\begin{proof}
Let $G$ be a maximal planar graph on $n=13$ vertices. Then from the property that every edge contributes one 4-cycle, then $C_4(G)\geq 3n-6=33$. If there is a separating $4$-cycle, then we have $C_4(G)\geq 34.$ Moreover, it can be checked that if 
$G$ contains a separating $3$-cycle, then $G$ contains a separating $4$-cycle as well, and hence the bound still holds. Indeed, let $T$ be a separating 3-cycle with vertices $x_1,x_2$ and $x_3$. Then one of the two regions contains at least two vertices. Without loss of generality, assume that the interior of the triangle contains two vertices and let $z\in N(x_1)\cap N(x_2)$ which is in the interior of the triangle and $(x_1,\ z,\ x_2, \ x_1)$ is a face. Now consider the $4$-cycle defined by $(x_1,\ z,\ x_2, \ x_3,\ x_1)$ which is clearly a separating $4$-cycle. 

Observe that the maximal planar graph shown in Figure~\ref{2}(right), which is also mentioned in~\cite{hakimi}, contains $34$ $4$-cycles and this makes the bound attainable. Now to complete proof of the claim, it is enough to verify that every $13$-vertex maximal planar graph contains a separating $4$-cycle.

We prove that if $\Delta(G)\geq 6$, then $G$ contains a separating $4$-cycle. Indeed, suppose for contradiction $G$ has no separating $4$-cycle. Let $v\in V(G)$ such that $d(v)\geq 6$. Denote $d(v)=k$ and $(v_1,\ v_2,\ v_3,\ \dots,\ v_k,\ v_1)$ be the unique $k$-cycle induced by vertices in $N(v)$. Denote the cycle by $\mathcal{C}$. If $\mathcal{C}$ is chordal, then $G$ contains a separating $3$-cycle, and consequently has a separating $4$-cycle which is a contradiction. 

Since $G$ is a maximal planar graph, each edge in $G$ is incident to two distinct $3$-faces. Hence for $v_iv_{i+1}\in E(\mathcal{C})$ there is a vertex $v\neq w_i\in V(G)$ such that $(v_i,\ w_i,\ v_{i+1},\ v_i)$ forms a $3$-face. We may assume that $w_i$ is unique to $v_iv_{i+1}$. To see why, let  $v_iv_{i+1}, \ v_jv_{j+1}\in E(\mathbf{C})$ distinct edges such that $w_i=w_j$. Without loss of generality assume $v_i\neq v_j$. In this case, we have a $4$-cycle $(v_i,\ w_i,\ v_{j+1},\ v,\ v_i)$ which separates the vertices $v_j$ and $v_{j+2}$. But this is a contradiction. 

In fact, $G$ does not exist if $k\geq 7$. Otherwise, we need more than $13$ vertices and this is not possible as $G$ is a $13$-vertex maximal planar graph. Assume that $k=6$, i.e., $d_G(v)=6$. In this case, $d_G(v_i)=5$ for each $i\in[6]$, and hence $w_iw_{i+1}$, for $i\in [5]$. We may assume that $d_G(w_1)>4$. If $w_1\in N(w_3)$, then $d(w_2)=4$, and we have a separating $4$-cycle, the cycle induced by $N(w_2)$. If $w_1\in N(w_4)$, then it can be checked that either $d(w_2)=4$ or $d(w_3)=4$, which results in a contradiction.    
\end{proof}

\end{document}